\nopagenumbers

\headline={\hss\tenrm\folio\hss}
\magnification
\magstep1
\baselineskip 18pt
\hsize=160truemm
\vsize=245truemm

\parindent 0.0cm
\par \centerline {{\bf SOME RESULTS ON THE JANOWSKI'S STARLIKE }}
\par \centerline  {{\bf FUNCTIONS OF COMPLEX ORDER }}
\vskip 0.3cm
\par \centerline {By : YASAR POLATOGLU AND METIN BOLCAL }
\vskip 0.3cm
\par Department of Mathematics of Faculty of Science and Arts ˜stanbul Kltr University
\vskip 0.3cm
\parindent 0.5cm
\par {\bf ABSTRACT : } The purpose of this paper is to give exact value of the radius
of starlikeness and distortion theorem, Koebe domain for the class of Janowski's
starlike functions of complex order. We note that the class of Janowski's starlike
functions of complex order contain many interesting subclasses of univalent functions.
\vskip 0.3cm
\par {\bf INTRODUCTION :} Let  $\Omega $  be the family of functions $\omega (z)$  regular in the unit disc $ D=\{ z \Big \vert  \vert z \vert < 1 \} $ and satisfying the conditions
$\omega (0)=0 \>\>,\>\>\vert \omega (z) \vert < 1 $ for $z \in D $.
\parindent 0.7cm
\par Next for arbitrary fixed numbers $A$ and $B$ , $-1 < A \leq 1 \> $,$\> -1 \leq B < A \> $, denote by
$P(A,B)$ the family of functions
$$p(z) = 1+p_{1}z+p_{2}z^{2}+\cdots \leqno(1.1)$$
\vskip 0.3cm
\parindent 0.0cm
regular in $D$ such that $p(z)$ in $P(A,B)$ if and only if
$$ p(z) = { {1+A\omega (z)} \over {1+B\omega (z)}} $$
\vskip 0.3cm
for some functions $\omega (z) \in \Omega $ and every $z \in D $. (This class was
introduced by Janowski [9] ).
\parindent 0.5cm
\par  Moreover let $S^{*}(A,B,b)\>\>$$(b\not = 0,complex )$ denote the family of functions.
$$ f(z) =z+a_{2}z^{2}+a_{3}z^{2}+\cdots \leqno (1.2) $$
\vskip 0.3cm
\parindent 0.0cm
regular in $D$ and such that $f(z) \in S^{*}(A,B,b) $ if and only if
$$ \Bigg [ 1+{1 \over b } (z {{f^{'}(z)} \over {f(z)} }-1) \Bigg ]=p(z)\quad ,\quad (b\not =0 \>,\> complex ) \leqno (1.3) $$
\vskip 0.3cm
for some $p(z)$ in $P(A,B)$ and all $z$ in $D$.
\parindent 0.5cm
\par Finally, we use $p$ denote the class of functions
$$p_{1}(z) = 1+c_{1}z+c_{2}z^{2}+\cdots \leqno(1.4)$$
\parindent 0.0cm
\vskip 0.3cm
which are analytic in $D$ and have a positive real part in $D$.
\parindent 0.5cm
\par It is to be note that special selections $A,B$ and $b$ lead to familiar sets of univalent
functions. Therefore the sets of univalent functions are listed by the following
\parindent 0.0cm
\par {\bf 1) }$S^{*}(1,-1,1) $ is the class of starlike functions (well known class ) [1,2,6].
\par {\bf 2) }$S^{*}(1,-1,b) $ is the class of starlike functions of complex order.Introduced by Wiatrowski [7].
\par {\bf 3) }$S^{*}(1,-1,1-\beta)\>, \> 0\leq \beta < 1 ,\> $ is the class of starlike functions of order $\beta $.This class was introduced by Robertson [5].
\par {\bf 4) }$S^{*}(1,-1,e^{-i\lambda } Cos\lambda )\>, \> \vert \lambda \vert < {\pi \over 2 } $ is the class of  $\lambda $-spirallike functions introduced by Spacek [3].
\par {\bf 5) }$S^{*}(1,-1,(1-\beta )e^{-i\lambda } Cos\lambda) \>, \> 0\leq \beta < 1 \>,\> \vert \lambda \vert < {\pi \over 2 }$
is the class of $\lambda $-spirallike functions of order $\beta $.This class was introduced by Libera [8].
\parindent 0.7cm
\vskip 0.3cm
\par The expression
$\Bigg [ 1+{1 \over b } (z {{f^{'}(z)} \over {f(z)} }-1) \Bigg ]$ is denoted by $ST(b)$,then
\vskip 0.3cm
\parindent 0.0cm
\par {\bf 6) }$S^{*}(1,0,b) $ is the set defined by $ \vert ST(b)-1 \vert <1 $.
\par {\bf 7) }$S^{*}(\beta ,0,b) $ is the set defined by $ \vert ST(b)-1 \vert <\beta \>,\> 0 \leq \beta <1 $.
\par {\bf 8) }$S^{*}(\beta ,-\beta ,b) $ is the set defined by $ \Bigg\vert {{ST(b)-1}\over {ST(b)+1} } \Bigg\vert <\beta \>,\> 0 \leq \beta <1 $.
\par {\bf 9) }$S^{*}(1,(-1+{1 \over M }),b) $ is the set defined by $ \vert ST(b)-M \vert <M $.
\par {\bf 10) }$S^{*}(1-2\beta ,-1,b) $ is the set defined by $ Re ST(b)>\beta \>,\> 0\leq \beta <1 $.
\parindent 0.5cm
\vskip 0.3cm
\par {\bf II. THE RADIUS OF STARLIKENESS FOR THE CLASS $S^{*}(A,B,b)$ }
\vskip 0.3cm
\parindent 0.7cm
\par From the definition of the class of $S^{*}(A,B,b)$ we easily obtain the following lemma.
\vskip 0.3cm
\parindent 0.5cm
\par {\bf LEMMA 2.1.} Let $f(z) \in S^{*}(A,B,b)$,then
$$\Bigg \vert z { {f^{'}(z) }\over {f(z)}}-{
{1-\Big [ B^{2}-b(AB-B^{2})\Big ]r^{2} } \over {1-B^{2}r^{2} } } \Bigg \vert \leq
{ {\vert b \vert (A-B)r } \over {1-B^{2}r^{2} }}. $$
\vskip 0.3cm
\par {\bf Proof: } Let $p(z) \in P(A,B)$ then
$$\Big \vert p(z)-{{1-ABr^{2}} \over {1-B^{2}r^{2} }} \Big \vert  \leq
{ { (A-B)r } \over {1-B^{2}r^{2} }}.\leqno(2.1) $$
\parindent 0.0cm
\vskip 0.3cm
The relation (2.1) was proved by Janowski [9]. Therefore from the definition of the class
$S^{*}(A,B,b)$ we can write
$$ \Bigg \vert \Bigg [ 1+{1 \over b } (z {{f^{'}(z)} \over {f(z)} }-1) \Bigg ]-
{{1-ABr^{2}} \over {1-B^{2}r^{2} }} \Bigg \vert \leq
{ { (A-B)r } \over {1-B^{2}r^{2} }} .\leqno (2.2) $$
\vskip 0.3cm
After the simple calculations from the relation (2.2) we obtain the desired result of the lemma.
\parindent 0.5cm
\par {\bf THEOREM 2.1.} The radius of starlikeness of the class $S^{*}(A,B,b)$ is
$$r_{s}= { 2 \over {\vert b \vert (A-B)+ \sqrt {{\vert b \vert}^{2} (A-B)^{2} +4 [B^{2}
+(AB-B^{2})Re \> b ] }}}. $$
\parindent 0.0cm
This radius is sharp. Because the extremal function is


$$f_{*}(z)= \cases{
$z(1+Bz)$^{{b(A-B)} \over B}$$  &$B \not= 0$  \cr
$e$^{bAz}$$ &B=0 \cr }$$

$$z={
{r \Big ( r-{{\overline b } \over b }\Big )^{1/2} } \over {1-r({{\overline b } \over b })^{1/2} }}$$
\parindent 0.5cm
\vskip 0.3cm
\par {\bf Proof: } From the Lemma 2.1. we have
\vskip 0.3cm
$$ Re \> z {{f^{'}(z)} \over {f(z)} }\geq {{1- \vert b \vert (A-B)r-[B^{2}+(AB-B^{2})]r^{2} } \over {1-B^{2}r^{2} }}.\leqno(2.3) $$
\vskip 0.3cm
\parindent 0.7cm
\par Hence for $r<r_{s}$ the first side of the preceding inequality is positive,this implies that
$$r_{s}= { 2 \over {\vert b \vert (A-B)+ \sqrt {{\vert b \vert}^{2} (A-B)^{2} +4 [B^{2}
+(AB-B^{2})Re \> b ] }}}. $$
\par Also note that the inequality (2.3) become an equality for the function$f_{*}(z)$.
\vskip 0.3cm
\par {\bf I. }For $A=1\quad ,\quad B=-1,$
$$ r_{s}={1 \over {\vert b \vert +\sqrt {{\vert b \vert}^{2}-2 Re\>b +1}}}.$$
\vskip 0.3cm
\parindent 0.0cm
This is the radius of starlikeness for the class of starlike functions of complex order was obtained by M.A.Nasr and M.K.Aouf [4].
\parindent 0.0cm
\par {\bf In this case;}
\par If we give the special values to $\quad b\quad $we obtain the radius of starlikeness for the corresponding classes. These radius had been found the authers [1],[3],[5],[8],[9].
\parindent 0.7cm
\par {\bf II.}For $A=1\quad ,\quad B=0,$
$$r_{s}={1 \over {\vert b \vert }}.$$
\parindent 0.0cm
\vskip 0.3cm
\par {\bf  In this case;}
\par Under the conditions $\vert \lambda \vert < {\pi \over 2 }\>\>,\>\> 0\leq\alpha <1 .$
\parindent 0.0 cm
\par {\bf (i)} $b=1 \quad r_{s} =1$ \quad ,\quad {\bf (ii)}$b=1-\alpha \quad r_{s}={1 \over {1-\alpha }}$\quad ,\quad {\bf (iii)}$ b=e^{-i\lambda }.Cos \lambda \quad  r_{s} ={1\over {Cos \lambda }} \quad ,$
\vskip 0.3cm
\par {\bf (iv)}$b=(1-\alpha) e^{-i\lambda }.Cos \lambda \quad  r_{s}={1 \over { (1-\alpha)Cos \lambda }}$
\parindent 0.7cm
\vskip 0.3cm
\par {\bf III.} For $A=\beta\quad ,\quad B=0\quad , \quad (0 \leq \beta <1 ). $
$$r_{s}= {1 \over {\beta \vert b \vert }}.$$
\parindent 0.0cm
\vskip 0.3cm
\par {\bf In this case; }
\par Under the conditions $\vert \lambda \vert < {\pi \over 2 }\quad ,\quad 0 \leq \alpha < 1 .$
\parindent 0.0cm

\par {\bf (i)} $b=1 \quad r_{s} ={1\over \beta } \quad ,\quad ${\bf (ii)}$b=(1-\alpha) \quad r_{s}={1 \over {\beta (1-\alpha) }}\quad ,$
\vskip 0.3cm
\par {\bf (iii)}$ b=e^{-i\lambda }.Cos \lambda \quad  r_{s} ={1\over {\beta Cos \lambda }} \quad ,\quad $ {\bf (iv)}$b=(1-\alpha) e^{-i\lambda }.Cos \lambda \quad  r_{s}={1 \over {\beta (1-\alpha)Cos \lambda }}$
\vskip 0.3cm
\parindent 0.7cm
\par {\bf IV.} For $A=\beta\quad ,\quad B=-\beta  \quad , \quad (0 \leq \beta <1 ). $
$$ r_{s}={1 \over {\beta [\vert b \vert +\sqrt {{\vert b \vert}^{2}-2Re \> b+1}]}}.$$
\parindent 0.0cm
\vskip 0.3cm
\par {\bf In this case ;}
\parindent 0.0cm
\par {\bf (i)} $b=1 \quad r_{s} ={1\over \beta } \quad ,\quad ${\bf (ii)}$b=(1-\alpha) \quad r_{s}={1 \over {\beta (1+2\alpha) }}\quad ,$
\vskip 0.3cm
\par {\bf (iii)}$ b=e^{-i\lambda }.Cos \lambda \quad  r_{s} ={1\over {\beta [Cos \lambda +\vert Sin \lambda \vert ] }} \quad , $
\vskip 0.3cm
\par {\bf (iv)}$b=(1-\alpha) e^{-i\lambda }.Cos \lambda \quad  r_{s}={1 \over {\beta [(1-\alpha)Cos \lambda + \sqrt {1-(1-{\alpha }^{2})Cos \lambda } ] }}.$
\vskip 0.3cm
\parindent 0.7cm
\par {\bf V.} For $A=1-2\beta \quad , \quad B=-1 .$
$$ r_{s}={1 \over {(1-\beta) \vert b \vert +\sqrt {1+2(1-\beta )Re\> b+{\vert b \vert}^{2}\beta ^{2}}}}.$$
\parindent 0.0cm
\vskip 0.3cm
\par {\bf In this case ; }
\parindent 0.0cm
\par {\bf (i)} $b=1 \quad r_{s} ={1\over {(1-\beta)+\sqrt {{\beta }^{2}-2\beta +3 } } }\quad ,\quad ${\bf (ii)}$b=(1-\alpha) \quad r_{s}={1 \over {2(1-\alpha)(1-\beta )+1 }}\quad ,$
\vskip 0.3cm
\par {\bf (iii)}$ b=e^{-i\lambda }.Cos \lambda \quad  r_{s} ={1\over {(1-\beta ) Cos \lambda + \sqrt {1+({\beta }^{2}+2\beta +2)Cos^{2} \lambda}}} \quad , $
\vskip 0.3cm
\par {\bf (iv)}$b=(1-\alpha) e^{-i\lambda }.Cos \lambda \quad  r_{s}={1 \over {(1-\beta)(1-\alpha )Cos^{2} \lambda + \sqrt {1+[(1-\alpha)^{2}(1-\beta )^{2}+2(1-\alpha)(1-\beta)]Cos^{2} \lambda }}}.$
\vskip 0.3cm
\parindent 0.7cm
\par {\bf VI.} For $A=-1,B=({1 \over M}-1) .$
$$ r_{s}={1 \over {\vert b \vert(2-{1 \over M})+\sqrt {{\vert b \vert }^{2}{(2-{1 \over M})}^{2}+4{1 \over M}({1 \over M}-1)Re\> b +1}}}$$
\parindent 0.0cm
\vskip 0.3cm
\par {\bf In this case ; }
\parindent 0.0cm
\par {\bf (i)} $b=1 \quad r_{s} ={1\over {(2-{1 \over M})^{2}+\sqrt {5{1 \over{ M^{2}} }  - 4{1 \over M}+1}}}$
\vskip 0.3cm
\par {\bf (ii)}$b=(1-\alpha) \quad r_{s}={1 \over {(1-\alpha)(2-{1 \over M} )+\sqrt { (1-\alpha)^{2}(2-{1 \over M} )^{2} +4{1 \over M}({1 \over M}-1)(1-\alpha )+1}}}$
\vskip 0.3cm
\par {\bf (iii)}$ b=e^{-i\lambda }.Cos \lambda \quad  r_{s}={1 \over {(2-{1 \over M} )Cos \lambda +\sqrt { ({4 \over {M^{2}}}-{8 \over M}+4)Cos^{2} \lambda + 1}}}$
\vskip 0.3cm
\par {\bf (iv)}$b=(1-\alpha) e^{-i\lambda }.Cos \lambda \quad
 r_{s}={1 \over {(2-{1 \over M} )(1-\alpha )Cos \lambda +\sqrt { ({4 \over {M^{2}}}-{8 \over M}+4)(1-\alpha )^{2}Cos^{2} \lambda + 1}}}$
\vskip 0.3cm
\parindent 0.5cm
\par {\bf III. THE ESTIMATION OF $\vert f(z) \vert $ IN $S^{*}(A,B,b)$}
\vskip 0.3cm
\parindent 0.7cm
\par In this section we shall give the estimation of $\vert f(z) \vert $ and the Koebe domain for the class of
$S^{*}(A,B,-b)$.
\par {\bf THEOREM 3.1.} If $f(z) \in S^{*}(A,B,b)$,then
\parindent 0.0cm
$$F(r;-A,-B,\vert b \vert ) \leq \vert f(z) \vert \leq F(r;A,B,\vert b \vert ) \leqno(3.1) $$
\vskip 0.3cm
where
$$F(r;A,B,\vert b \vert )=\cases{
$r(1+Br)$^{{{\vert b \vert (A-B)}\over B}}$$ &if $B\not =0$ ; \cr
$re$^{\vert b \vert Ar}$$ &if $B=0$ \cr }$$


\vskip 0.3cm
This bound are sharp. Because the extremal function is

$$f_{*}(z)=\cases {
$z(1+Bz)$^{{b(A-B)} \over B }$$ &if$B\not = 0 $ ; \cr
$ze$^{bAz}$$ &if$B=0$. \cr}$$

\parindent 0.7cm
\vskip 0.3cm
\par {\bf Proof :} Since  $f(z) \in S^{*}(A,B,b)$ we have
$$1+{1 \over b }(z{{f^{'}(z)} \over {f(z)}}-1)=p(z)\quad ,\quad p(z) \in P(A,B) \leqno (3.2)$$
\vskip 0.3cm
\parindent 0.0cm
and simple calculations from the equality (3.1) we obtain.
$$ f(z)=zExp \Bigg [{ \int_{0}^{z}} { {b(p(\xi)-1)}\over z} d\xi \Bigg ]  $$
\vskip 0.3cm
Therefore
$$\vert  f(z) \vert =\vert z \vert Exp \Bigg [Re \Big({ \int_{0}^{z}} { {b(p(\xi)-1)}\over {\xi }} d\xi \Big )\Bigg ] \leqno (3.3) $$
\vskip 0.3cm
substituting $\xi =zt$, we obtain
$$\vert  f(z) \vert =\vert z \vert Exp \Bigg [Re \Big({ \int_{0}^{1}} { {b(p(zt)-1)}\over t} dt \Big )\Bigg ]. \leqno (3.4) $$
\vskip 0.3cm
On the other hand from lemma 2.1 it follows that
$$\max\limits_{\vert zt \vert =rt } \Bigg ({ {b(p(zt)-1)}\over t} \Bigg ) =
{ {\vert b \vert (A-B)r} \over {1+Brt}};\leqno (3.5)$$
\vskip 0.3cm
then after integration we obtain the upper bounds in (3.1) similarly obtain the bounds
on the left-hand side of (3.1),which shows that the proof of the theorem is complete.
\parindent 0.7cm
\vskip 0.3cm
\par {\bf (I)}For $A=1 \quad ,\quad B=-1$
$${r \over {(1+r)^{2\vert b \vert }}}\leq \vert f(z) \vert \leq {r \over {(1-r)^{2\vert b \vert }}}.$$
\vskip 0.3cm
\par This is the distortion for the class of starlike functions of complex order.
\vskip 0.3cm
\parindent 0.0cm
\par {\bf In this case; }
\par {\bf (i) } For $b=1$;
$${r \over {(1+r)^{2}}}\leq \vert f(z) \vert \leq {r \over {(1-r)^{2}}}.$$
\vskip 0.3cm
this is the distortion for the class of starlike functions $\quad $(This is well known result,
\par A.W.Goodman, univalent functions vol1,page 140 )
\vskip 0.3cm
\par {\bf (ii) } For $b=(1-\alpha )\quad ,\quad 0\leq \alpha <1$;
$${r \over {(1+r)^{2(1-\alpha)}}}\leq \vert f(z) \vert \leq {r \over {(1-r)^{2(1-\alpha )}}}.$$
\vskip 0.3cm
This result is the distortion for the class of starlike function of order $\alpha $. This result was obtained by M.S.Robertson [5].
\vskip 0.3cm
\par {\bf (iii) } For $b=e^{-i\lambda }Cos \lambda \quad ,\quad \vert \lambda \vert <{\pi \over 2 }$;
$${r \over {(1+r)^{2Cos \lambda }}}\leq \vert f(z) \vert \leq {r \over {(1-r)^{2Cos \lambda }}}.$$
\vskip 0.3cm
This is the distortion for the class of $\lambda $-spirallike functions.
\vskip 0.3cm
\par {\bf (iv) } For $b=(1-\alpha )e^{-i\lambda }Cos \lambda \quad ,\quad \vert \lambda \vert <{\pi \over 2 }\quad , \quad 0 \leq \alpha <1 $;
$${r \over {(1+r)^{2(1-\alpha )Cos \lambda }}}\leq \vert f(z) \vert \leq {r \over {(1-r)^{2(1-\alpha )Cos \lambda }}}.$$
\vskip 0.3cm
This is the distortion for the class of $\lambda $-spirallike functions of order $\alpha$.
\vskip 0.3cm
\parindent 0.7cm
\par {\bf (II)}For $A=\beta \quad ,\quad B=-\beta $
$${r \over {(1+\beta r)^{2\vert b \vert }}}\leq \vert f(z) \vert \leq {r \over {(1-\beta r)^{2\vert b \vert }}}.$$
\vskip 0.3cm
\parindent 0.0cm
\par {\bf In this case; }
\par {\bf (i) } For $b=1$;
$${r \over {(1+\beta r)^{2}}}\leq \vert f(z) \vert \leq {r \over {(1-\beta r)^{2}}}.$$
\vskip 0.3cm
\par {\bf (ii) } For $b=(1-\alpha )\quad ,\quad 0\leq \alpha <1$;
$${r \over {(1+\beta r)^{2(1-\alpha)}}}\leq \vert f(z) \vert \leq {r \over {(1-\beta r)^{2(1-\alpha )}}}.$$
\vskip 0.3cm
\par {\bf (iii) } For $b=e^{-i\lambda } Cos \lambda \quad ,\quad \vert \lambda \vert <{\pi \over 2 }$;
$${r \over {(1+\beta r)^{2Cos \lambda }}}\leq \vert f(z) \vert \leq {r \over {(1-\beta r)^{2Cos \lambda }}}.$$
\vskip 0.3cm
\par {\bf (iv) } For $b=(1-\alpha )e^{-i\lambda} Cos \lambda \quad ,\quad \vert \lambda \vert <{\pi \over 2 }\quad , \quad 0 \leq \alpha <1 $;
$${r \over {(1+\beta r)^{2(1-\alpha )Cos \lambda }}}\leq \vert f(z) \vert \leq {r \over {(1-\beta r)^{2(1-\alpha )Cos \lambda }}}.$$
\vskip 0.3cm
\parindent 0.7cm
\par {\bf (III)}For $A=1 \quad ,\quad B=0 $
$${r \over {e^{\vert b \vert r}}}\leq \vert f(z) \vert \leq re^{\vert b \vert r} $$
\vskip 0.3cm
\parindent 0.0cm
\par {\bf In this case; }
\vskip 0.3 cm
\par {\bf (i) } For $b=1$;
$${r \over {e^{ r}}}\leq \vert f(z) \vert \leq re^{ r} $$
\vskip 0.3cm
\par {\bf (ii) } For $b=(1-\alpha )$;
$${r \over {e^{(1-\alpha) r}}}\leq \vert f(z) \vert \leq re^{(1-\alpha) r} $$
\vskip 0.3cm
\par {\bf (iii) } For $b=e^{-i\lambda }Cos \lambda \quad ,\quad \vert \lambda \vert < {\pi \over 2 } $;
$${r \over {e^{ rCos \lambda }}}\leq \vert f(z) \vert \leq re^{ rCos \lambda } $$
\vskip 0.3cm
\par {\bf (iv) } For $b=(1-\alpha )e^{-i\lambda }Cos \lambda \quad ,\quad 0 \leq \alpha < 1 \quad ,\quad \vert \lambda \vert < {\pi \over 2 } $;
$${r \over {e^{ (1-\alpha )rCos \lambda }}}\leq \vert f(z) \vert \leq re^{(1-\alpha ) rCos \lambda } $$
\vskip 0.3cm
\parindent 0.7cm
\par {\bf (IV)}For $A=\beta  \quad ,\quad B=0 $
$${r {e^{-\vert b \vert \beta r}}}\leq \vert f(z) \vert \leq re^{\vert b \vert \beta r} $$
\vskip 0.3cm
\parindent 0.0cm
\par {\bf In this case; }
\vskip 0.3 cm
\par {\bf (i) } For $b=1$;
$$r {e^{- \beta r}}\leq \vert f(z) \vert \leq re^{ \beta r} $$
\vskip 0.3cm
\par {\bf (ii) } For $b=(1-\alpha )\quad ,\quad 0 \leq \alpha < 1 $;
$$r e^{-(1-\alpha) \beta r}\leq \vert f(z) \vert \leq re^{(1-\alpha) \beta r} $$
\vskip 0.3cm
\par {\bf (iii) } For $b=e^{-i\lambda }Cos \lambda \quad ,\quad \vert \lambda \vert < {\pi \over 2 } $;
$$r  e^{ -(\beta Cos \lambda )r }\leq \vert f(z) \vert \leq re^{ (\beta Cos \lambda )r } $$
\vskip 0.3cm
\par {\bf (iv) } For $b=(1-\alpha )e^{-i\lambda }Cos \lambda \quad ,\quad 0 \leq \alpha < 1 \quad ,\quad \vert \lambda \vert < {\pi \over 2 } $;
$$r e^{ -(\beta (1-\alpha )Cos \lambda )r }\leq \vert f(z) \vert \leq re^{(\beta (1-\alpha ) Cos \lambda )r } $$
\vskip 0.3cm
\parindent 0.7cm
\par {\bf (V)}For $A=1-2\beta  \quad ,\quad B=-1$
$${r \over {(1+r)^{2\vert b \vert (1-\beta )}}}\leq \vert f(z) \vert \leq {r \over {(1-r)^{2\vert b \vert (1-\beta )}}}.$$
\vskip 0.3cm
\parindent 0.0cm
\par {\bf In this case; }
\par {\bf (i) } For $b=1$;
$${r \over {(1+r)^{2(1-\beta )}}}\leq \vert f(z) \vert \leq {r \over {(1-r)^{2(1-\beta )}}}.$$
\vskip 0.3cm
\par {\bf (ii) } For $b=(1-\alpha )$;
$${r \over {(1+r)^{2(1-\alpha)(1-\beta )}}}\leq \vert f(z) \vert \leq {r \over {(1-r)^{2(1-\alpha )(1-\beta )}}}.$$
\vskip 0.3cm
\par {\bf (iii) } For $b=e^{-i\lambda }Cos \lambda $;
$${r \over {(1+r)^{2(1-\beta )Cos \lambda }}}\leq \vert f(z) \vert \leq {r \over {(1-r)^{2(1-\beta )Cos \lambda }}}.$$
\vskip 0.3cm
\par {\bf (iv) } For $b=(1-\alpha )e^{-i\lambda }Cos \lambda  $;
$${r \over {(1+r)^{2(1-\alpha )(1-\beta )Cos \lambda }}}\leq \vert f(z) \vert \leq {r \over {(1-r)^{2(1-\alpha )(1-\beta )Cos \lambda }}}.$$
\vskip 0.3cm
\vskip 0.3cm
\parindent 0.7cm
\par {\bf (VI)}For $A=1  \quad ,\quad B={1 \over M }-1$
$${r \over {\Big [1+(1-{1\over M })r\Big ]^{\vert b \vert (2-{1 \over M } )}}}\leq \vert f(z) \vert \leq {r \over {\Big [1-(1-{1\over M })r \Big ]^{\vert b \vert (2-{1\over M } )}}}.$$
\vskip 0.3cm
\parindent 0.0cm
\par {\bf In this case; }
\vskip 0.3cm
\par {\bf (i) } For $b=1$;
$${r \over {\Big [1+(1-{1\over M })r\Big ]^{ (2-{1 \over M } )}}}\leq \vert f(z) \vert \leq {r \over {\Big [1-(1-{1\over M })r \Big ]^{ (2-{1\over M } )}}}.$$
\vskip 0.3cm
\par {\bf (ii) } For $b=(1-\alpha ) $;
$${r \over {\Big [1+(1-{1\over M })r\Big ]^{(1-\alpha ) (2-{1 \over M } )}}}\leq \vert f(z) \vert \leq {r \over {\Big [1-(1-{1\over M })r \Big ]^{(1-\alpha ) (2-{1\over M } )}}}.$$
\vskip 0.3cm
\par {\bf (iii) } For $b=e^{-i\lambda }Cos \lambda $;
$${r \over {\Big [1+(1-{1\over M })r\Big ]^{ (2-{1 \over M } )Cos \lambda }}}\leq \vert f(z) \vert \leq {r \over {\Big [1-(1-{1\over M })r \Big ]^{ (2-{1\over M } )Cos \lambda }}}.$$
\vskip 0.3cm
\par {\bf (iv) } For $b=(1-\alpha )e^{-i\lambda }Cos \lambda $ ;
$${r \over {\Big [1+(1-{1\over M })r\Big ]^{ (2-{1 \over M } )(1-\alpha ) Cos \lambda }}}\leq \vert f(z) \vert \leq {r \over {\Big [1+(1-{1\over M })r \Big ]^{ (2-{1\over M } )(1-\alpha )Cos \lambda }}}.$$
\vskip 0.3cm
\parindent 0.5cm
\par {\bf IV. KOEBE DOMAIN FOR THE CLASS {\bf $S^{*}(A,B,b)$ } }
\parindent 0.7 cm
\par In this section we shall give the Koebe domain for the class $S^{*}(A,B,b)$
under the condition of the definition of Koebe domain.[see 1. page 113-114 ].
\par From the theorem 3.1. we have;
$$ \vert f(z) \vert \geq r(1-Br)^{ {{\vert b \vert (A-B)} \over B}} \quad \quad for \quad B\not = 0 $$
$$ \vert f(z) \vert \geq re^{-\vert b \vert Ar } \quad \quad for \quad B \not = 0 $$
\vskip 0.3cm
\parindent 0.0 cm
then from the definition of Koebe domain we obtain
$$ R= \lim_{r\rightarrow 1^{-} } [r(1-Br)^{ {{\vert b \vert (A-B)} \over B}}]=(1-B)^{ {{\vert b \vert (A-B)} \over B}} \quad \quad for \quad B\not = 0 \leqno (4.1) $$
$$ R= \lim_{r\rightarrow 1^{-} } [ re^{-\vert b \vert Ar } ]= e^{-\vert b \vert A } \quad \quad for \quad B = 0  \leqno (4.2)$$
\vskip 0.3cm
\parindent 0.7cm
\par {\bf (I)}For $A=1 \quad ,\quad  B=-1 \quad \quad
R= {1 \over {4^{\vert b \vert }}}  $
\vskip 0.3cm
\parindent 0.0cm
\par {\bf In this case ;}
\vskip 0.3cm
\par {\bf (i)}$b=1 \quad \quad R={1 \over 4 } $. This is well known result (Koebe-${1 \over 4 }$theorem [1,page 115 ])
\vskip 0.3cm
\par {\bf (ii)}$b=1-\alpha  \quad \quad R={1 \over {4^{(1-\alpha )}} } $. This result was obtained by M.S.Robertson
\par [1. page 115,[5]].
\vskip 0.3cm
\par {\bf (iii)}$b=e^{-i\lambda }Cos \lambda  \quad \quad R={1 \over {4^{Cos \lambda }} } $. This is the Koebe domain for the class of $\lambda $-spirallike function.
\vskip 0.3cm
\par {\bf (iv)}$b=(1-\alpha )e^{-i\lambda }Cos \lambda  \quad \quad R={1 \over {4^{(1-\alpha )Cos \lambda }} } $. This is the Koebe domain for the class of $\lambda $-spirallike functions of order $\alpha $.
\vskip 0.3cm
\parindent 0.7cm
\par {\bf (II)}For $A=1 \quad ,\quad  B=0 \quad \quad
R= {1 \over {e^{\vert b \vert }}} . $
\vskip 0.3cm
\par {\bf (III)}For $A=\beta  \quad ,\quad  B=0 \quad \quad
R= {1 \over {e^{\vert b \vert \beta }}} . $
\vskip 0.3cm
\par {\bf (IV)}For $A=\beta  \quad ,\quad  B=-\beta \quad \quad
R= {1 \over {(1-\beta )^{2\vert b \vert \beta }}} . $
\vskip 0.3cm
\par {\bf (V)}For $A=1  \quad ,\quad  B=-1+{1 \over M }\quad \quad
R= {
1 \over
{\Big (2-{1 \over M} \Big )^{{\vert b \vert (2-{1 \over M })} \over{1-{1 \over M }}}}} .$
\vskip 0.3cm
\par {\bf (VI).}For $A=1-2\beta  \quad ,\quad  B=-1\quad \quad
R= {1 \over {4^{\vert b \vert (1-\beta ) }}} . $
\vskip 0.8cm
\parindent 0.0cm
\par \centerline {R E F E R E N C E S }
\vskip 0.5 cm
\par {\bf [1] } A.W.Goodman;Univalent functions Vol I and Vol II.Mariner publishing
\par $\quad $Company.Inc.Tampa Florida 1983.
\par {\bf [2] } CHR.Pommerenke;Univalent functions.Vandenhoeck,ruprecht in G"ttingen,1975 .
\par {\bf [3] } L.Spacek;Prispeek k teori funki prostych,Casopis pest Math.Fys.62 (1933) 12-19.
\par {\bf [4] } M.A.Nasr and M.K.Aouf;Starlike functions of complex order.Jour of Natural science
\par $\quad $ and Mathematics Vol 25.No 1 (1985) 1-12.
\par {\bf [5] } M.S.Robertson;On the theory of univalent functions.Ann.of Math. 37(1936)374-408
\par {\bf [6] } P.L.Duren;Univalent functions,Springer-Verlag 1983.
\par {\bf [7] } P.Wiatrowski;The Coeffecient of a certain family of holomorphic functions.Zeszyty.
\par $\quad  $ Nauk.Math.przyord.ser II. Zeszyt.(39)Math (1971) 57-85.
\par {\bf [8] } R.J.Libera;Some radius of convexity problem.Duke Math. J.31(1964) 143-157.
\par {\bf [9] } W.Janowski;Extremal problem for a family of functions with positive real part and
\par $\quad $ for some related families.Ann.Polon.Math 23(1970) 159-177.

\end